\documentclass[12pt]{amsart}

\usepackage{amsthm}
\usepackage{amsmath}
\usepackage{amsfonts}
\usepackage{amssymb} 
\usepackage{amscd}
\usepackage{hyperref}
\newtheorem*{theorem*}{Theorem}
\newtheorem*{lemma*}{Lemma}
\theoremstyle{definition}

\begin{document}

\title{Another Riemann Hypothesis Equivalent}

\author{Richard Pell}

 \email{richard.pell@frontier.com}

\begin{abstract}Identities involving M\"obius function values $(\mu(j), \mu(k))$ are used to 
generate a Riemann Hypothesis equivalent. 
\end{abstract}

\keywords{Riemann Hypothesis, M\"obius function, Mertens function, Dirichlet series}
\subjclass[2010]{11M26, 11M41}

\maketitle

The purpose of this article is to present a Riemann Hypothesis equivalent
with such a simple structure that anyone with just a modest background in
number theory can study it.  We use a linear combination of identities involving
the M\"obius function to show the following.    The key advantage is that the expression 
inside the parentheses assumes only the values zero and one.

\begin{theorem*} Let $\mu(n)$ be the M\"obius function and let $[\phantom{x}]$ denote the greatest integer function. The Riemann Hypothesis is equivalent to the statement
$$
\overline{\lim}_{n\to\infty} \frac{\log\left|\sum_{j=1}^n\sum_{k=1}^n \left(\left[\frac{n^2}{jk}\right] - 2\left[\frac{n^2}{2jk}\right]\right)\mu(k)\mu(j)\right|}{\log n} \le 1.
$$
\end{theorem*}

\begin{proof} We start with the following.

\begin{lemma*} Let $j, k, m$ be positive integers. Then
\begin{enumerate}
\item $$
\left[\frac{\left[\frac{m}{j}\right]}{k}\right] = \left[\frac{m}{jk}\right].$$
\item
$$
\sum_{k=1}^m \left[\frac{m}{k}\right] \mu(k) = 1.
$$
\item
$$
\sum_{k=1}^m \left[\frac{m}{jk}\right] \mu(k) = 1
$$
when $j\le m$.
\item 
$$
\sum_{k=1}^m\left( \left[\frac{m}{jk}\right] - 2\left[\frac{m}{2jk}\right]\right) \mu(k)= \begin{cases} -1 \text{ if } j\le m/2 \\ 1 \text{ if } m/2 < j \le m. \end{cases}
$$
\end{enumerate}
\end{lemma*}
\begin{proof}
\begin{enumerate}
\item See \cite{Card} (p. 163, Lemma 6).
\item Meissel's identity \cite{Meis}.
\item Replace $m$ with $[m/j]$ in (2). Then use (1).
\item This follows from (3).
\end{enumerate}
\end{proof}

Part 4 of the lemma implies that
\begin{align*}
-\sum_{j=1}^n  \mu(j)
&=\sum_{j=1}^n\sum_{k=1}^{n^2} \left(\left[\frac{n^2}{jk}\right] - 2\left[\frac{n^2}{2jk}\right]\right)\mu(k)\mu(j)\\
& = \sum_{k=1}^{n^2}\sum_{j=1}^{n} \left(\left[\frac{n^2}{jk}\right] - 2\left[\frac{n^2}{2jk}\right]\right)\mu(k)\mu(j)\\
&= \sum_{k=1}^{n}\sum_{j=1}^{n} \left(\left[\frac{n^2}{jk}\right] - 2\left[\frac{n^2}{2jk}\right]\right)\mu(k)\mu(j)\\
&\phantom{aaaa}+  \sum_{k=n+1}^{[n^2/2]}\sum_{j=1}^{n} \left(\left[\frac{n^2}{jk}\right] - 2\left[\frac{n^2}{2jk}\right]\right)\mu(k)\mu(j)\\
& \phantom{aaaa}+ \sum_{k=[n^2/2]+1}^{[n^2]}\sum_{j=1}^{n} \left(\left[\frac{n^2}{jk}\right] - 2\left[\frac{n^2}{2jk}\right]\right)\mu(k)\mu(j)\\
& = \sum_{k=1}^{n}\sum_{j=1}^{n} \left(\left[\frac{n^2}{jk}\right] - 2\left[\frac{n^2}{2jk}\right]\right)\mu(k)\mu(j)\\
&\phantom{aaaa}+  \sum_{k=n+1}^{[n^2/2]}(-1)\mu(k)  \qquad \text{ (from Lemma, part 4)}\\
&\phantom{aaaa}+  \sum_{k=[n^2/2]+1}^{[n^2]}(1)\mu(k)   \qquad \text{ (from Lemma, part 4)}.
\end{align*}
Let $M(x) = \sum_{k\le x} \mu(k)$ be the Mertens function. Subtract $\sum_{k=1}^n \mu(k)$ from both sides of the above and rearrange to obtain
\begin{align*}
M\left(n^2\right) - 2 M([n^2/2]) & = -\sum_{j=1}^n\sum_{k=1}^n \left(\left[\frac{n^2}{jk}\right] - 2\left[\frac{n^2}{2jk}\right]\right) \mu(j)\mu(k) - 2\sum_{k=1}^n \mu(k)\\
&= -\sum_{j=1}^n\sum_{k=1}^n \left(\left[\frac{n^2}{jk}\right] - 2\left[\frac{n^2}{2jk}\right]\right) \mu(j)\mu(k) + O(n).
\end{align*}

$M(n^2) - 2M([n^2/2])$ is the $n^2$th partial sum of the coefficients of
$$
(1-2^{1-s}) \sum_{k\ge 1} \mu(k) k^{-s}.
$$
The convergence of this series for $Re(s)>1/2$ implies that $\zeta(s)$ has no zeros for $Re(s) > 1/2$, which is the Riemann Hypothesis.
(Note that all of the zeros of $1-2^{1-s}$ lie on the vertical
line whose real part equals $1$. On that same line, the reciprocal
of the zeta function has no poles (PNT).)

Conversely, if the Riemann Hypothesis holds, then 
$$
M(x) = O_{\epsilon}(x^{1/2+\epsilon}).
$$
Therefore, 
$$
M(x) - 2M([x/2]) =  O_{\epsilon}(x^{1/2+\epsilon}),
$$
which implies that the series converges for $Re(s)>1/2$.

The series cannot converge for $Re(s) < 1/2$ because $\zeta(s)$ has zeros with $Re(s) = 1/2$. 

It follows that the Riemann Hypothesis is equivalent to the statement that
$$
\overline{\lim}_{x\to\infty} \frac{\log|M(x) - 2 M([x/2]|}{\log (x)} \le 1/2.
$$

Therefore,
$$
\overline{\lim}_{n\to\infty} \frac{\log|M(n^2) - 2 M([n^2/2]|}{\log (n^2)} \le 1/2
$$
is implied by the Riemann Hypothesis. 
It remains to show that using the subsequence of perfect squares is sufficient to imply the Riemann Hypothesis. 

Let $\epsilon > 0$ and suppose 
$$
|M(n^2) - 2M([n^2/2])| \le Cn^{1+\epsilon}$$
 for all $n\ge N_0$.  Let $x\ge N_0^2$ be arbitrary and choose $n\ge N_0$ with $n^2 \le x < (n+1)^2$. 
Then 
$$
|M(x) - M(n^2)| \le 2n \quad \text{ and } \quad| M([x/2]) - M([n^2/2]) | \le n.
$$
Therefore, 
$$
|(M(x)-2M([x/2])) - (M(n^2)-2M([n^2/2]))| \le 4n,
$$
so
$$
M(x)-2M([x/2]) = O(x^{1/2+\epsilon}).
$$
This implies that
$$
\overline{\lim}_{x\to\infty} \frac{\log|M(x) - 2 M([x/2])|}{\log x} \le 1/2,
$$
which implies the Riemann Hypothesis, as stated above.
\end{proof}

\section*{Acknowledgements}  Thanks to Jeffrey Lagarias for help when I first began
studying the Mertens function.
Special thanks to Lawrence Washington for helping
with the preparation of this manuscript.


\begin{thebibliography}{99}
\bibitem{Card} J.-P. Cardinal, ``Symmetric matrices related to the Mertens function,'' {\it Linear Algebra Appl.} 432 (2010), pp. 161-172.
\bibitem{Meis} E. Meissel, ``Observationes quaedam in theoria numerorum,'' {\it J. reine angew. Math.} 48 (1854), pp. 301-316.
\end{thebibliography}
\end{document}